\documentclass[12pt]{article}
\usepackage{amsfonts}
\usepackage{amsthm}
 \usepackage{amsmath}
\author{M.J. Heath}
\title{A note on a construction of J.F. Feinstein.}

\begin{document}
\maketitle
\begin{abstract}
In \cite{F}  J.F. Feinstein constructed a compact plane set $X$ such that $R(X)$ has no non-zero, bounded point derivations but is not
weakly amenable. In the same paper he gave an example of a separable uniform algebra $A$ such that every point in the character space of  $A$ is a peak
point  but $ A$ is not weakly amenable. We show that it is possible to modify the construction in order to produce examples which are also
regular.  
\end{abstract}
\baselineskip=18pt
\section{Introduction}
\newtheorem{Theorem}{Theorem}[section]

\newtheorem{Corollary}[Theorem]{Corollary}
\newtheorem{Example}[Theorem]{Example}
   
\newtheorem{Observation}[Theorem]{Observation}
\newtheorem{Lemma}[Theorem]{Lemma}
\newtheorem{Proposition}[Theorem]{Proposition}
\newtheorem{claim}{Claim}
\theoremstyle{definition}
\newtheorem{Definition}[Theorem]{Definition}
\newtheorem{Claim}[claim]{Claim}
\def\sskip{\vskip 0.3 cm\noindent}
A uniform algebra $A$ on  a compact (Hausdorff) space $X$ is said to be \emph{regular on $X$} if for any point $x$ in $X$ and any
compact subset $K$ of $X\setminus\{x\}$ there is a function $f$
in $A$ such that $f(x)=1$ and $f$ is zero on $K$. We call $A$ \emph{regular} if it is regular on its character space. $A$ is said to be
\emph{trivial} if it is $C(X)$ - the uniform algebra of all continuous functions on $X$. The first example of a non-trivial
regular uniform algebra was given by McKissick \cite{Mc} (see also \cite{K}; \cite[chapter 37]{S}): the example was $R(X)$ for a
compact plane set $X$.
\paragraph{}The notion of weak amenability was introduced in \cite{B}. A commutative Banach algebra $A$ is said to be \emph{weakly
amenable} if there are no non-zero,
continuous derivations from $A$ into the any commutative Banach $A$-bimodule. It is proved in \cite {B} that this is equivalent to there
being no non-zero continuous derivations into the dual module $A^\prime$. 
\paragraph{}As point derivations may be regarded as derivations into 1-dimensional, commutative Banach modules it is a necessary condition for weak
amenability that there be no non-zero, bounded point derivations. However, this condition is not sufficient, even for uniform algebras:
in \cite{F} J.F. Feinstein constructed a compact plane set $X$ such that the uniform algebra $R(X)$ has no non-zero, bounded point
derivations but $R(X)$ is not weakly amenable. It was not clear whether such an example could also be regular. In this note we 
show that it can.
\paragraph{Notation} Throughout this paper $Q$ will refer to the compact plane set  $\{x+iy:x,y\in[-1,1]\}$.
\paragraph{} For a plane set $X$ and a function $f\in C(X)$, $\vert f\vert_X$ will be the uniform norm of $f$ on $X$, $\sup\{\vert
f(z)\vert : z\in X\}$. For a compact plane set $X$ we denote by $R_0(x)$ the set of restrictions to $X$ of rational functions with poles
off $X$. Hence the uniform algebra $R(X)$ is the uniform closure of $R_0(X)$ in $C(X)$.
\paragraph{}If $D$ is a disc in the plane then $r(D)$ shall refer to its radius. 

\paragraph{} Let $\mu$ be a complex measure on a compact plane set $X$ such that the bilinear functional on $R_0(X)\times R_0(X)$ defined by
\begin{displaymath}(f,g)\mapsto \int_xf^\prime(x)g(x)d\mu(x)
\end{displaymath} is bounded. Then, as in \cite{F}, we may extend by continuity to $R(X)\times R(X)$ and obtain a continuous derivation
$D$ from $R(X)$ to $R^\prime(X)$ such that for any $f$ and $g$ in $R_0(X)$ we have,
\begin{displaymath}D(f)(g)=\int_Xf^\prime(x)g(x)d\mu(x).
\end{displaymath}
\paragraph {}In the next section, we strengthen K\"orner's \cite{K} version of McKissick's Lemma to allow greater control over the centres and radii of
the discs removed. Using this we modify Feinstein's \cite{F} construction so that $R(X)$ is regular. In fact we prove the following theorem. 

\begin{Theorem} \label{thm}For each $C>0$ there is a compact plane set $X$ obtained by deleting from $Q$ a countable union of Jordan domains such that
$\partial Q$ is a subset of $X$, $R(X)$ is regular and has no non-zero, bounded point derivations and, for all $f, g$ in $R_0(X),$
\begin{displaymath} \bigg\vert\int_{\partial Q}f^\prime(z)g(z)\textrm{d}z\bigg\vert\le C\vert f\vert_X\vert g\vert_X.
\end{displaymath}
\end{Theorem}
 If we let $X$ be a compact plane set constructed as in Theorem \ref{thm} we have, by the discussion above, a non-zero continuous
 derivation $D$ from $R(X)$ to $R^\prime(X)$ such that
 \begin{displaymath}D(f)(g)=\int_{\partial Q}f^\prime(z)g(z)dz
 \end{displaymath}
 for $f,g\in R_0(X)$. So $R(X)$ is not weakly amenable. 
  
  \section{The Construction}
  Our main new tool will be the following theorem, which is a variant on McKissick's result in \cite{Mc} (see also \cite{K}).
  
  \begin{Theorem}\label{Mck} For any $C_0 >0$ there is a compact plane set $X_1$ obtained by deleting from $Q$ a countable union of open
  discs $(D_n)_{n=1}^\infty$ whose closures are in $\textrm{\emph {int}} (Q)$ such
  that $R(X_1)$ is regular and, letting $s_n$ be the distance of $D_n$ to $\partial Q$,
  \begin{displaymath}\sum_{n=1}^\infty \frac{r(D_n)}{s_n^{\phantom n 2}}<C_0.
  \end{displaymath}
  \end{Theorem}
  In order to prove Theorem \ref{Mck} we require a series of lemmas which are variants on the results of \cite{K}
  
  The following is  \cite[Lemma 2.1]{K}.
  
  \begin{Lemma}\label{Mc3}If $N\ge2$ is an integer and $h_N(z)=1/(1-z^N)$ then the following hold:
  \begin{itemize}
  \item[(i)]$\vert h_N(z)\vert\le2\vert z\vert^{-N}$ for $\vert z\vert^N\ge2$;
  \item[(ii)]$\vert1-h_N(z)\vert\le2\vert z\vert^N$ for $\vert z\vert^N\le2^{-1}$;
  \item[(iii)]$h_N(z)\ne0$ for all $z$.
  \end{itemize}
  Further if $(8\log N)^{-1}>\delta>0$ then;
  \begin{itemize}
  \item[(iv)]$\vert h_N(z)\vert\le2\delta^{-1}$ provided only that $\vert z-w\vert\ge\delta N^{-1}$ whenever $w^N=1$.
  \end{itemize}
  \end{Lemma}
  The following is a variant of  \cite[Lemma 2.2]{K}.
  \begin{Lemma}\label{Mc4} If in Lemma \ref{Mc3} we set $N=n2^{2n}$ with $n$ sufficiently large then:
  \begin{itemize}
  \item[(i)]$\vert h_N(z)\vert\le(n+1)^{-4}$ for $\vert z\vert\ge1+2^{-(2n+1)}$;
  \item[(ii)]$\vert1-h_N(z)\vert\le(n+1)^{-4}$ for $\vert z\vert \le1-2^{-(2n+1) }$;
  \item[(iii)]$h_N\ne0$ for all $z$;
  \item[(iv)]$\vert h_N(z)\vert\le n^{-4}2^{2n+1}$ provided only that $\vert z-w\vert \ge n^{-5}2^{-4n}$ whenever $w^N=1$.
   \end{itemize}
   \end{Lemma}
\begin{proof}Parts (i), (ii) and (iii) are the corresponding parts of \cite{K}, 2.2. Part (iv) follows on putting $\delta=n^{-3}2^{-2n}$
\end{proof}
From this we obtain the following variant on \cite[Lemma 2.3]{K} 
\begin{Lemma} \label{Mc5}Provided only that n is sufficiently large we can find a finite collection $A(n)$ of disjoint open discs and a rational function
$g_n$ such that, letting $s_0(\Delta)=\textrm{dist}(\Delta ,\mathbb{R}\cup i\mathbb{R})$, for a disc $\Delta$, the following hold: 
\begin{itemize}
\item[(i)] $\sum_{\Delta\in A(n)}\frac {r(\Delta)}{s_0(\Delta)^2}< n^{-2}$ and so $\sum_{\Delta\in A(n)}r(\Delta)\le n^{-2}$;
\item[(ii)]the poles of $g_n$ lie in $\bigcup_{\Delta\in A(n)}\Delta$;
\item[(iii)] $\vert g_n(z)\vert\le(n+1)^{-4}$ for $\vert z\vert\ge1-2^{-(2n+1)}$;
\item[(iv)] $\vert1-h_N(z)\vert\le(n+1)^{-4}$ for $\vert z\vert \le1-2^{-(2n-1) }$;
\item[(v)] $\vert g_n(z)\vert\le n^{4}2^{2n+1}$ for $z\not\in\bigcup_{\Delta\in A(n)}\Delta$;
\item[(vi)]$g_n(z)\ne0$ for all z;
\item[(vii)]$\bigcup _{\Delta\in A(n)}\Delta\subset\big\{z:\-2^{-(2n-1)}\le\vert z\vert\le1-2^{-(2n+1)}\big\}$.
\end{itemize}
\end{Lemma}
\begin{proof}Let $N=n2^{2n}$, $\omega=\exp(2\pi/N)$, $g_n=h_N(\omega^{-\frac{1}{2}}(1-2^{-2n})^{-1}z)$. If we take $A(n)$ to be the
collection of discs with radii $n^{-5}2^{-4n}$ and centres $(1-2^{-2n})\omega^{r+\frac{1}{2}}, [0\le r\le N-1],$ results (ii)-(vii) are
either trivial or follow from Lemma \ref{Mc4} by on scaling by a factor of $\omega^\frac{1}{2}(1-2^{-2n})$.
To show part (i), consider first those discs with centres $(1-2^{-2n})\omega^{r+\frac{1}{2}}, (0\le r\le \frac{N-1}{8})$   \\
For such a disc $\Delta$ we have:
\begin{eqnarray*}s_0(\Delta)&=&(1-2^{-2n})\sin\Big(\frac{(r+\frac{1}{2})\pi}{n2^{2n}}\Big)-n^{-5}2^{-4n}\\
&\ge&\frac{(r+\frac{1}{2})\pi}{n2^{2n+2}}-n^{-5}2^{-4n}\\
&\ge& \Big(\frac{1}{2}\Big)\frac{2r+1}{n2^{2n}}
\end{eqnarray*}.
So
\begin{displaymath} \frac{r(\Delta)}{s_0(\Delta)^2} 
\le n^{-5}2^{-4n}   \frac{4n^2 2^{4n}}{(2r+1)^2}
=\frac{4n^{-3}}{(2r+1)^2}
\end{displaymath}
and
\begin{displaymath}\sum_{k=1}^{\frac{N}{8}-1}\frac{r(\Delta_k)}{s_0(\Delta_k)}\le 4n^{-3}\sum_{k=0}^{\frac{N}{8}-1}(2k+1)^{-2}
\le 4n^{-3}\sum_{k=0}^ \infty k^{-2}\le\frac{n^{-2}}{8}
\end{displaymath}
provided only that $n> K=32\sum_{r=0}^\infty r^{-2}$.
So, by symmetry,
\begin{displaymath}\sum_{\Delta\in A(n)}\frac{r(\Delta)}{s_0(\Delta)^2}\le n^{-2}
\end{displaymath} provided only that $n$ is sufficiently large.
\end{proof}
Multiplying the $g_n$ together as in \cite{K}, we obtain the following.
\begin{Lemma}\label{Mc6}Given any $\varepsilon>0$ there exists an $m=m(\varepsilon)$ such that if (adopting the notation of Lemma \ref{Mc5}) we let
$f_n=(m!)^{-4}\prod_{r=m}^ng_r$ and $\{\Delta_k\}$ be a sequence enumerating the discs of $\bigcup_{r=m}^\infty A(r)$ then the following hold:
 \begin{itemize}
 \item[(a)]  $\sum_{k=1}^\infty\frac{r(\Delta_k)}{s_0(\Delta_k)^2}<\varepsilon$ and so then $\sum_{k=1}^\infty r(\Delta_k) <\varepsilon$;
 \item[(b)] The poles of the $f_n$ lie in $\bigcup_{k=1}^\infty\Delta_k$;
 \item[(c)] The sequence $\{f_n\}$ tends uniformly to zero on $\{z\in \mathbb{C}:\vert z\vert\ge1\}$.
 \end{itemize}
\end{Lemma}

\begin{proof}
Observe that
 \begin{displaymath}\sum_{k=1}^\infty r(\Delta_k)\le \sum_{k=1}^\infty \frac{r(\Delta_k)}{s_0(\Delta_k)^2}=\sum_{l=m}^\infty\sum_{\Delta\in
 A(r)}\frac{r(\Delta)}{s_0(\Delta)^2}\le\sum_{l=m}^\infty l^{-2}<\varepsilon
\end{displaymath}
provided only that $m(\varepsilon)>2\varepsilon^{-1}+1$. Thus conclusions (a) and (b) are easy to verify. To prove (c), set $K=\prod_{r=1}^\infty
(1+(r+1)^{-4})$ and observe that, provided $m(\varepsilon)$ is large enough that  $(m+1)!^4>2(m+1)^6 2^{2(m+1)}+2(m+1)^2$, then if we let
$z\not\in\bigcup_{k=1}^\infty\Delta_k$,  a simple induction gives
\begin{eqnarray*}\vert f_n(z)\vert&\le&(n+1)!^{-4}\textrm{ for }1-2^{-(2n+1)}<\vert z\vert,\\
\vert f_n(z)\vert&\le& n^{-2}\le\prod_{r=m}^n (1+(r+1)^{-4})\le K \textrm{ for }1-2^{-(2n-1)}<\vert z\vert<1-2^{-(2n+1)},\\
\vert f_n(z)\vert&\le&\prod_{r=m}^n(1+(r+1)^{-4})\le K\textrm{ for }\vert z \vert< 1-2^{-(2n-1)}.\\
\end{eqnarray*}
  Using the trivial equality
  \begin{displaymath}\vert f_{n+1}(z)-f_n(z)\vert=\vert f_n(z)\vert\vert1-g_{n+1}(z)\vert,
  \end{displaymath}
we see that for $z\not\in\bigcup_{k=1}^\infty\Delta_k$,
\begin{eqnarray*}\vert f_{n+1}(z)-f_n(z)\vert&\le& K(n+1)^{-4}\textrm{ for }\vert z\vert\le1-2^{-(2n+1)},\\
\vert f_{n+1}(z)-f_n(z)\vert&\le&(n+1)!^{-4}(1+(n+1)^4 2^{2n+1})\le(n+1)^{-2},\\
\textrm{ for }1-2^{-(2n+1)}<\vert z\vert
\end{eqnarray*}

Thus $ \vert f_{n+1}(z)-f_n(z)\vert\le K(n+1)^{-2}$ for all $z\not\in\bigcup_{k=1}^\infty\Delta_k$ and, by for example the Weierstrass $M$ test, $f_n$
converges uniformly to $f$ say.
To see that $f(z)\ne0$ for $\vert z\vert<1$, $z\not\in\bigcup_{k=1}^\infty\Delta_k$, note that if $\vert z\vert\le1-2^{(2n-1)}$ then $f_n(z)\ne0$, and
\begin{displaymath} \sum_{r=n+1}^\infty\vert1-g_r(z)\vert\le\sum_{r-n+1}^\infty(r+1)^{-4}<\infty.
\end{displaymath}
So by a basic result on infinite products (see, for example, \cite[15.5]{Ru} ),
\begin{displaymath}f(z)=f_n(z)\prod_{r=n+1}^\infty g_r(z)\ne0.
\end{displaymath}
\end{proof}
Whence by dilation and translation we obtain the following.
 \begin{Lemma} \label{Mc2}Given any closed disc $D$, with centre $a$ and radius $r$, and any $\varepsilon>0$ we can find a sequence of open discs
$\{\Delta_k\}$ and a sequence of rational functions $\{f_n\}$ such that, letting $s_1(\Delta):= \textrm{dist}(\Delta_k, a+\mathbb{R}\cup i\mathbb{R})$ for a disc
$\Delta$, the following hold:
 \begin{itemize}
 \item[(a)] $\sum_{k=1}^\infty\frac{r(\Delta_k)}{s_1(\Delta_k)^2} <\varepsilon$ and so $\sum_{k=1}^\infty R(\Delta_k)<\varepsilon$;
 \item[(b)] The poles of the $f_n$ lie in $\bigcup_{k=1}^\infty\Delta_k$;
 \item[(c)] The sequence $\{f_n\}$ tends uniformly to zero on $(\mathbb{C}\setminus D)\setminus\bigcup_{k=1}^\infty\Delta_k$.
 \end{itemize}
 \end{Lemma}
 We are now ready to prove Theorem \ref{Mck}. 
\paragraph{Proof of Theorem \ref{Mck}.}
Let $\{D_l\}_{l=1}^\infty$ be an enumeration of all closed discs of centre $z$ and radius $r$  with $z\in \mathbb Q +i\mathbb Q$ and $r\in \mathbb Q^+$ such
that, letting $K= \{-1-i, -1+i, 1-i, 1+i\}$, one of the following holds:
\begin{eqnarray} z&\in& \textrm{int}(Q)\textrm{ and } r<\textrm{dist}(z,\partial Y);\\
z&\in& \partial Q\setminus K \textrm{ and } r<\textrm{dist}(z,K);\\
z&\in& K \textrm{ and } r<1.
\end{eqnarray}
We apply Lemma \ref{Mc2} with $D = D_l$ and $\varepsilon = \varepsilon_l$ where
\begin{eqnarray*} \varepsilon_l&<&2^{-l-1}C_0\textrm{dist}(D_l,\partial Q)^2\textrm{ if }D_l \textrm{ is of type }(1),\\
\varepsilon_l&<&2^{-l-1}C_0\textrm{dist}(D_l,K)^2\textrm{ if }D_l \textrm{ is of type }(2),\\
\varepsilon_l&<&2^{-l-1}C_0\textrm{ if }D_l \textrm{ is of type }(3),
\end{eqnarray*}
to obtain $(\Delta_ {l,n})$ and $(f_{l,k})$. Let $\{U_N\}_{N=1}^\infty$ be a sequence enumerating the $\Delta_{n,a}$ and
\begin{displaymath}X_0=Q\setminus \bigcup_{N=1}^\infty U_N.
\end{displaymath}
We have
\begin{displaymath} \sum_{N=1}^\infty \frac{r(U_N)}{s(U_N)^2}=\sum_{l=1}^\infty\sum_{n=1}^\infty \frac
{r(\Delta_{l,n})}{s(\Delta_{l,n})^2}<\sum_{l=1}^\infty C_02^{-l}=C_0.
\end{displaymath}
Given any point $z$ in $X_0$ and any compact set $B\subset X_0$ there exists $D_l$ with $z\in D_l$ and $B\cap D_l = \emptyset$. Hence
$f_l:=\lim_{k\rightarrow\infty}f_{l,k}\in R(X_0)$ has $f_l(z)\ne 0$ and $f_l(B)\subset\{0\}$, so $R(X_0)$ is regular.
\begin{flushright}$\square$\end{flushright} 

\paragraph{} In order to prove Theorem \ref{thm} we need some further lemmas. The first of these lemmas is trivial.
\begin{Lemma}\label{triv} Let $X$, $Y$ be compact plane sets with $X\subset Y$. If $R(Y)$ is regular, then the same is true for $X$.
\end{Lemma}
The next two results are essentially the same as those used in  \cite{F}.
\begin{Lemma}\label{int} Let $(\chi_n)$ be a sequence of Jordan domains whose closures are contained in $Q$. Set
$X_2=Q\setminus\bigcup_{n=1}^\infty \chi_n$. Let $s_n$ be the distance from $\chi_n$ to $\partial Q$ and let $c_n$ be the length of the boundary of
$\chi_n$. Let $f$ and $g$ be in $R_0(X_2)$. Then
\begin{displaymath}\bigg\vert\int_{\partial Q}f^\prime(z)g(z)\textrm d z\bigg\vert\le2\vert f\vert _X\vert g\vert_X\sum_{n=1}^\infty
\frac{c_n}{s_n^2}.
\end{displaymath}
\end{Lemma}
\begin{proof}The argument of \cite[2.1]{F} applies.
\end{proof}
\begin{Lemma}\label{seq}
Let $X$ be a compact subset of $Q$. Suppose that there is a sequence of real numbers $L_n\in(0,1)$ such that 
$L_n\rightarrow 1$ and, for each $n$ $R(X\cap L_nQ)$ has no non-zero bounded point derivations. Then $R(X)$ has no non-zero bounded 
point derivations.
\end{Lemma}
\begin{proof} The argument of \cite[2.2]{F} applies.
\end{proof}
The following is \cite[2.3]{F} .
\begin{Lemma}\label{sub} Let $X$, $Y$ be compact plane sets with $X\subset Y$. If $R(Y)$ has no non-zero, bounded point derivations, then the same
is true for $X$.
\end{Lemma}
The following result was proved by Wermer \cite{W}.
\begin{Proposition}\label{prop} Let $D$ be a closed disc in $\mathbb{C}$ and let $\varepsilon>0$. Then there is a sequence of open discs $(U_k)\subset D$
such that $R(D\setminus \bigcup _{k=1}^\infty U_k)$ has no non-zero bounded point derivations but such that the sum of the radii of the
discs $U_k$ is less than $\varepsilon$.
\end{Proposition}
\begin{Corollary}\label{cor}Let $Y$ be a square set of the form $rQ+z$  and $\delta>0$. Then there is a sequence of Jordan domains
$(\chi_l)\subset Y$
such that $R(Y\setminus \bigcup _{k=1}^\infty \chi_l)$ has no non-zero bounded point derivations but such that the sum of the
lengths of the boundaries  of the  $\chi_l$ is less than $\delta$.
\end{Corollary}
\begin{proof}
Apply the previous proposition to any closed disc containing $Y$ with $\varepsilon < \frac{\delta}{2\pi}$. Then, by Lemma \ref{sub}, 
letting $(\chi_l)$ be a sequence enumerating all non-empty sets of the form $U_k\cap \textrm{int}(Y)$ will suffice.
\end{proof}
\paragraph{Proof of Theorem \ref{thm}}
Let $C>0$. Set $L_n=n/(n+1)$. Applying Corollary \ref{cor} to $L_nQ$  we may choose Jordan domains $\chi_{n,k}\subset L_nQ$ such that
$R(L_nQ\setminus\bigcup _{k=1}^\infty\chi_{n,k})$ has no non-zero bounded point derivations and such that the sum of the lengths of the boundaries of the 
$(\chi_{n,k})_{k=0}^\infty$ is less than $2^{-(n+1)}C(1-L_n)^2$.\\
Set 
\begin{displaymath}X_2=Q\setminus\bigcup_{n,k}\chi_{n,k}
\end{displaymath}
and let
\begin{displaymath}X_1=Q\setminus\bigcup_nD_n
\end{displaymath}
be the result of applying Theorem \ref{Mck} with $C_0=C/(4\pi)$. Finally set $X=X_1\cap X_2$ Then Lemma \ref{triv} gives that $R(X)$ is
regular and Lemmas \ref{seq} and \ref{sub} give that $R(X)$ has no non-zero bounded point derivations. Enumerating the  sets $\chi_{n,k}$ and $D_{n}$ as
$C_1,C_2, \dots$ we may apply Lemma \ref{int} to obtain the required estimate on the integral.  

\begin{flushright}$\square$
\end{flushright}
\section{Regularity and peak points.}
A point $x$ in the character space $X$ of a uniform algebra $A$ is said to be a \emph{peak point for} $A$ if there is $f\in A$ such that
$f(x)=1$ and $\vert f(y)\vert<1$ for all $y\in X\setminus\{x\}$ and $x$ is a \emph{point of continuity for} $A$ if, for every compact set $K\in X\{x\}$,
there is a function $f$ in $A$ such $f(x)=1$ and $f(K)\subset\{0\}$. Notice that a uniform algebra is regular if and only if every point in its character
space is a point of continuity. The following result, regarding systems of Cole root extensions (see \cite{C}), is \cite[2.8]{F2}
\begin{Proposition} Let $A$ and $B$ be uniform algebras such that $B$ is the result of applying a system of Cole root extensions to $A$. If $A$ is
regular so is $B$.

\end{Proposition}

\paragraph{}In \cite{F} an example of a  separable uniform algebra, $A$, such that every point of the character space is a peak point for $A$ and $A$ is
not  weakly amenable is obtained by first modifying the original example so that, every point, except possibly those on the outer boundary circle, is a
point of continuity, and then applying an appropriate system of Cole root extensions. Following essentially the same argument and noting the above
proposition we obtain the following.
 
\begin{Theorem} \label{pp} There exists a regular uniform algebra $A$ whose character space is metrizable such that every point of the character
space of $A$ is a peak point but such that $A$ is not weakly amenable. 
\end{Theorem} 
  
\paragraph{}We note that both this algebra and the algebra constructed in Theorem \ref{thm} have dense invertible group (by results in \cite{DF}) .

\paragraph{}We finish by noting that we do not know whether or not either of the uniform algebras we have constructed are strongly regular -  i.e. if, 
for each point $x$ in the character space,  the algebra of functions constant on a neighbourhood of $x$ is dense in the original algebra.

\sskip

{\sf  School of Mathematical Sciences

 University of Nottingham

 Nottingham NG7 2RD, England

 email: pmxmjh@nottingham.ac.uk}
\sskip
2000 Mathematics Subject Classification: 46J10, 46H20

\end{document}